\documentclass[12pt]{amsart}
\usepackage{amssymb,latexsym,amsthm,amsmath}

\theoremstyle{plain}

\newtheorem*{CF}{The Chung Feller Theorem}
\newtheorem{lemma}{Lemma}

\newtheorem*{prop}{Proposition}
\theoremstyle{remark}

\theoremstyle{definition}

\newcommand{\tri}{\triangleleft}
\begin{document}
\title[Chung-Feller Theorem]{An Investigation of\\ the Chung-Feller
Theorem}
\author{Eli~A. Wolfhagen}
\begin{abstract}
In this paper, we shall prove the Chung-Feller Theorem in several
ways. We provide an inductive proof, bijective proof, and proofs
using generating functions, and the Cycle Lemma of Dvoretzky and
Motzkin \cite{DM47}.
\end{abstract}
\maketitle
\section{Introduction}\label{S:intro}
The main focus of this paper is to prove the following result of
Chung and Feller \cite{CF49}:
\begin{CF}
The number of paths from $(0,0)$ to $(n,n)$, with steps $(0,1)$ and
$(1,0)$ and exactly $2k$ steps above the line $x=y$ is independent
of $k$, for every $k$ such that $0\leq k\leq n$. In fact, it is
equal to the $n$th Catalan number.
\end{CF}

\section{Preliminary Definitions and an Inductive Proof}\label{S:def}
A Dyck path is a path with steps $(0,1)$ and $(1,0)$ that starts at
the origin and ends at $(n,n)$ for some positive integer $n$. We can
study the number of paths indicated in the above Chung-Feller
theorem in a simpler form if we generalize Dyck paths as seen in
\cite{wF60}. A \emph{$k$-negative path} of length $2n$ is a path
from $(0,0)$ to $(2n,0)$, with steps $(1,1)$ and $(1,-1)$, such that
exactly $2k$ of these steps are below the horizontal axis.

A \emph{prime Dyck path} is a Dyck path that returns to the $x$-axis
only once at the end of the path. A \emph{negative prime Dyck path}
is a a prime Dyck path reflected about the $x$-axis. In general, any
\emph{negative} Dyck path is the reflection about the $x$-axis of a
Dyck path.

\begin{proof}[An Inductive Proof of the Chung-Feller Theorem]
As a base case, for $n=0$, there is only one path with $0$ steps. So
clearly there is only $C_0 = 1$ path with $0$ steps below the
$x$-axis.

Now, let $n > 0$ and assume that for all $i$ such that $0\leq i<n$,
the number of $l$-negative paths of length $2i$ is equal to $C_i$
for all $l\leq i$. Therefore, choose some $k \leq n$.

Any nonempty $k$-negative path either starts out with either a prime
Dyck path or a negative prime Dyck path. In the first case, the
$k$-negative path starts with an up-step followed by a Dyck path of
some length $2p-2$ which is then followed by a down-step and a
$k$-negative path of length $2n-2p$, for some $p\leq n-k$. On the
other hand, the second case deals with paths that start out with a
down-step, followed by a $(q-1)$-negative path of length $2q-2$
which is followed by an up-step and a $(k-q)$-negative path of
length $2n-2q$, for some $q\leq k$.

Let $\mathcal{N}_{\text{up}}$ be the number of $k$-negative paths of
length $2n$ that start with an up-step. Let
$\mathcal{N}_{\text{down}}$, similarly, be the number of
$k$-negative paths of length $2n$ that start with a down-step. Let
$\mathcal{N}$ be the total number of $k$-negative paths of length
$2n$.

If a path starts out with an up-step and a Dyck path of length
$2p-2$, then the number of such paths is the total number of Dyck
paths of length $2p-2$ multiplied by the number of $k$-negative
paths of length $2n-2p<2n$. By the inductive hypothesis, this means
that the number $N^{+}_p$ of such paths is $C_{p-1}C_{n-p}$. For the
total number of $k$-negative paths that start off with an up-step we
need to take the sum of $N^+_p$, for all $1\leq p\leq n-k.$
Therefore, we have that the total number of $k$-negative paths of
length $2n$ that start with an up-step is
\[\mathcal{N}_{\text{up}}=\sum_{p=1}^{n-k} N^+_p = \sum_{p=1}^{n-k}
C_{p-1}C_{n-p}.\]

Similarly, we define $N^-_q$ as the number of $k$-negative paths
that start out with a down-step and a negative Dyck path of length
$2q-2$. So $N^-_q$ is equal to the number of negative Dyck paths of
length $2q-2$ multiplied by the number of $(k-q)$-negative paths of
length $2n-2q<2n$. Obviously, since the negative Dyck paths are
reflections of (positive) Dyck paths, the number of negative Dyck
paths of length $2q-2$ is equal to the number of Dyck paths of
length $2q-2$. Therefore, we again have the sum of products of two
Catalan numbers, and the total number of $k$-negative paths of
length $2n$ that start with a down-step is
\[\mathcal{N}_{\text{down}} = \sum_{q=1}^k C_{q-1}C_{n-q} = \sum_{q=1}^k
C_{n-q}C_{q-1}.\]

Therefore the total number of $k$-negative paths of length $2n$ is
equal to
\begin{equation}\label{E:N} \mathcal{N}=\mathcal{N}_{\text{up}}+\mathcal{N}_{\text{down}} =
\sum_{p=1}^{n-k} C_{p-1}C_{n-p} + \sum_{q=1}^k
C_{n-q}C_{q-1};\end{equation} that is,
\[\mathcal{N} = (C_0C_{n-1} + C_1C_{n-2} + \cdots + C_{n-k-1}C_k)
+ (C_{n-1}C_0 + C_{n-2}C_1 + \cdots + C_{n-k}C_{k-1}).\] Reversing
the order of the summands in the second parentheses clearly gives
\begin{align*}\mathcal{N} &= (C_0C_{n-1} + \cdots +
C_{n-k-1}C_k)+(C_{n-k}C_{k-1} \cdots + C_{n-1}C_0)\\ &=
\sum_{i=0}^{n-1} C_iC_{n-i-1}.\end{align*}

Now, for any Dyck path of length $2n$ we can look at the first
nonempty prime path. It will have length $2i+2$ for some $i \geq 0$,
giving a total of $C_i$ such prime paths. Since the total length of
the entire path is $2n$, we must have $i \leq n-1$. Therefore, since
the rest of the path is an arbitrary Dyck path of length $2n-2i-2$
there are $C_{n-i-1}$ possible paths that can follow the initial
prime path of length $2i+2$. Therefore, the total number of Dyck
paths of length $2n$ that start with a prime path of length $2i+2$
is given by $C_iC_{n-i-1}$. Therefore if we sum over $0\leq i \leq
n-1$ we will get all possible Dyck paths of length $2n$, which is
given by $C_n$. Therefore $$ \mathcal{N} = \sum_{i=0}^{n-1}
C_iC_{n-i-1} = C_n,$$ and so the total number of $k$-negative paths
of length $2n$ is equal to $C_n$. \end{proof}

\section{A Bijective Proof}

Let $S_k$ denote the set of all $k$-negative paths of length $2n$.
Choose some $k$, such that $0 \leq k < n$.

For $s\in S_k$, find the last positive prime Dyck path $q$. Next
factor $s$ into $s=pqr$, where $p$ is the path up to the last
positive prime Dyck path and $r$ is the rest of the path after $q$.
Since $q$ is the last positive prime Dyck path, the remainder of the
path $r$ must be a negative Dyck path. This factorization is unique.

Like any other prime Dyck path, $q$ is composed of an up-step
followed by an arbitrary Dyck path $Q$ and a down-step. So $q$ can
be rewritten as $q=uQd$, where $u$ denotes an up-step and $d$
denotes a down-step. Define the function $\varphi_+: S_k \to
S_{k+1}$ by
\begin{equation}\label{E:phi} \varphi_+(s) = \varphi_+ (puQdr) = pdruQ, \;
\; \text{for any} s\in S_k.
\end{equation}

Given an $s\in S_k$, $\varphi_+(s) = pdruQ$ by \eqref{E:phi}. Since
$p$ begins and ends on the the $x$-axis, the step $d$ is below the
$x$-axis. As noted earlier $r$ is an arbitrary negative Dyck path,
so in the path $\varphi_+(s)$, $r$ starts and ends at height $-1$,
without going above it. Therefore the step $u$ is negative, starting
at height $-1$ and ends on the $x$-axis. Since $Q$ is a Dyck path it
contains no negative steps, so altogether $\varphi_+(s)$ has $2$
steps below the $x$-axis in addition to the number of negative steps
in path segments $p$ and $r$. Since the only negative steps in the
path $s\in S_k$ occur in $p$ and $r$, the number of negative steps
in path segments $p$ and $r$ is equal to $2k$. Thus $\varphi_+(s)$
has $2k+2$ steps below the $x$-axis and so is a member of $S_{k+1}$.

Now for any $\sigma \in S_{k+1}$ we can find the last negative prime
path $\omega$. Thus we can write $\sigma = \pi \omega \varrho$,
where $\pi$ is the path up to the last negative prime, and $\varrho$
is the rest of the path which is by construction positive. So since
$\omega$ is a negative prime path it can be written as $ d \Omega
u$, where $\Omega$ is an arbitrary negative Dyck path. Therefore any
$\sigma \in S_{k+1}$ can be uniquely written as $\sigma = \pi d
\Omega u \varrho$, where $\Omega$ is a negative path and $\varrho$
is a positive path. Therefore, if $s_1 = p_1uQ_1dr_1 \ne s_2 =
p_2uQ_2dr_2 \in S_k$ then either $p_1 \ne p_2$, $Q_1 \ne Q_2$ or
$r_1 \ne r_2$, so clearly since $\varphi_+(s_1) = p_1dr_1uQ_1,
\varphi_+(s_2)= p_2dr_2uQ_2$, since the factorization in $S_{k+1}$
is unique, $\varphi_+(s_1) \ne \varphi_+(s_2)$. Therefore,
$\varphi_+$ is injective.

Using the factorization on $S_{k+1}$ we can define the injection
$\varphi_- : S_{k+1} \to S_k$ which sends $\sigma = \pi d \Omega u
\varrho \mapsto \pi u \varrho d \Omega$. Since there are now two
fewer negative steps in $\varphi_-(\sigma)$ than in $\sigma$ itself,
due to the fact that the up-step now comes before the down-step,
$\varphi_-(\sigma) \in S_k$. Additionally, since there is an unique
decomposition in $S_k$, $\varphi_-$ is injective. Thus, by the
Scr\"{o}der-Bernstein theorem $|S_k| = |S_{k+1}|$. In fact,
$\varphi_- = \varphi_+^{-1}$, since $\varphi_-(pdruQ)=puQdr$ for all
$s = puQdr \in S_k$, so $\varphi_+$ is in fact a bijection.

Therefore, we can now bijectively prove the Chung-Feller theorem.
\begin{proof}[Bijective Proof] The cardinality of $S_k$ and
$S_{k+1}$ are equal, for all nonnegative $k\leq n-1$. Since there
are $n+1$ equal sets of paths of length $2n$ each with the same
cardinality, and the total number of paths of length $2n$ is
$\binom{2n}{n}$, we have the equation \[ (n+1)\mathcal{K} =
\binom{2n}{n},\] where $\mathcal{K}$ is the number of $k$-negative
paths of length $2n$ for any $k$ such that $0\leq k \leq n.$

Therefore $\mathcal{K} = \frac{1}{n+1} \binom{2n}{n} = C_n.$
\end{proof}

\section{A Generating Function Approach to the
Theorem}\label{S:genfn}

In general a generating function is a very useful tool used in
enumerative combinatorics. The generating function of an infinite
sequence $\{a_n\}$ can be thought of as the function for which the
coefficient of $x^n$ in the power series expansion about $x = 0$, is
$a_n$. That is, if $A(x)$ is the generating function for the
sequence $\{a_n\}$, then $A(x) = \sum_{n=0}^\infty a_nx^n$.

The generating function can also be obtained more formally by
weighting certain combinatorial objects by a variable $x$. For
instance, if $\{a_n\}$ represents the number of \emph{bloops} with
$n$ \emph{glurps}, then $A(x)$ is obtained by weighting each glurp
by $x$ and summing over all bloops. We will use both formulations in
this section to give a straight-forward proof of the Chung-Feller
theorem.

The generating function for the Catalan numbers is given by
\begin{equation}\label{E:cx} c(x) = \sum_{n=0}^\infty C_nx^n =
\frac{1-\sqrt{1-4x}}{2x},\end{equation} which is obtained by
weighting every step in a Dyck path by $\sqrt{x}$ and summing over
all possible Dyck paths.

Let us construct a generating function for arbitrary paths that end
on the $x$-axis by weighting each step by $\sqrt{x}$ and each step
below the $x$-axis by $\sqrt{tx}$, such that a path with $2n$ total
steps, $2k$ of which are below the $x$-axis is given the weight
$t^kx^n$. As described in section \ref{S:def}, the primes of a
$k$-negative path are either positive prime Dyck paths or negative
prime Dyck paths. Let $P_n$ denote the number of (positive) prime
Dyck paths of length $2n$.

Such a path of length $2n$ consists of an arbitrary Dyck path of
length $2n-2$ sandwiched between an up-step and down-step. So the
number of prime Dyck paths of length $2n$ is given by the $(n-1)$th
Catalan number, that is $P_n=C_{n-1}$. Therefore the generating
function for the prime Dyck paths is given by
\begin{equation}\label{E:p} p_+(x) =\sum_{n=1}^\infty P_nx^n =
\sum_{n=1}^\infty C_{n-1}x^n = x\sum_{n=0}^\infty C_nx^n =
xc(x).\end{equation}

Similarly, for negative prime Dyck paths the generating function,
now weighted by $\sqrt{xt}$ because each step in a negative prime
Dyck path is below the $x$-axis is given by
\begin{equation}\label{E:pn} p_-(x,t) = \sum_{n=1}^\infty P_n(tx)^n
= tx\sum_{n=0}^\infty C_n(tx)^n = txc(tx). \end{equation} So since
arbitrary paths can be factored into $l$ primes (either positive or
negative) for some $l \geq 0$, the generating function for such
paths is
\[N(x,t) = \sum_{l=0}^\infty (p_- + p_+)^l = \frac{1}{1-(p_- +
p_+)}.\] In this generating function the coefficient of $t^ix^j$ is
the number paths of total length $2j$ and with $2i$ steps below the
$x$-axis.

By the definition of $c(x)$ and equations \eqref{E:p} and
\eqref{E:pn} we get that
\begin{align*} N(t,x) &= \frac{1}{1-xc(x)-txc(tx)}\\
&= \frac{1}{\displaystyle{1-(x\frac{1-\sqrt{1-4x}}{2x} +
tx\frac{1-\sqrt{1-4tx}}{2tx})}}\\ &=
\frac{2}{\sqrt{1-4x}+\sqrt{1-4tx}}.
\end{align*} Rationalizing the denominator we see that
\begin{align}\label{E:Ngenfn}
N(t,x) &=\frac{\sqrt{1-4x}-\sqrt{1-4tx}}{2x(t-1)}\\
&=\frac{1-\sqrt{1-4x}-1+\sqrt{1-4xt}}{2x(1-t)}\notag \\
&=\frac{1}{1-t}\cdot \left( \frac{1-\sqrt{1-4x}}{2x} -
\frac{1-\sqrt{1-4xt}}{2x}\right)\notag \\ &=\frac{1}{1-t}\left(
c(x)-tc(tx)\right)\notag \\&= \frac{1}{1-t}\left( \sum_{n=0}^\infty
C_nx^n - \sum_{n=0}^\infty C_nt^{n+1}x^n\right)\notag\\
&=\sum_{n=0}^\infty \frac{1-t^{n+1}}{1-t} C_nx^n.\notag \end{align}

Using $N(t,x)$ and a simple algebraic identity we can prove the
Chung-Feller Theorem.

\begin{proof}
Because $$\frac{1-t^{n+1}}{1-t} = 1 + t + t^2 + \cdots + t^n,$$
equation \eqref{E:Ngenfn} can be rewritten as
\begin{equation}\label{E:Ntfn}
N(t,x)= \sum_{n=0}^\infty C_n \left(\sum_{k=0}^n t^k\right)x^n.
\end{equation}
Thus the coefficient of $t^kx^n$ in $N(t,x)$ is equal to $C_n$ for
$0\leq k\leq n.$ Therefore the number of $k$-negative paths of
length $2n$ is equal to $C_n$ for $0\leq k\leq n$.
\end{proof}

\section{A Proof by Reordering}
The Cycle Lemma of Dvoretzky and Motzkin \cite{DM47} is intricately
linked with our main theorem. In \cite{eW04}, the author provides
the following proof of the Cycle Lemma using an ordering which
naturally reveals the equidistribution of the Chung-Feller Theorem.

\begin{lemma}[Cycle Lemma] Given a sequence $\pi = a_1a_2a_3 \cdots
a_{n}$ of integers with $a_i \leq 1$ such that $\sum_{j=1}^n a_i = k
> 0$, there are exactly $k$ values of $i$ such that all partial sums
of the cyclic permutation $\pi_i = a_{i+1}a_{i+2}\cdots a_{2n+k}a_1
\cdots a_i$ are positive. \end{lemma}

This is a stronger claim than we need, so let us just restrict
ourselves to $a_i \in \{-1,1\}$. For ease of computation for any $p
\leq n$, let $s(p) = \sum_{j=1}^p a_i$. Given a sequence $\pi$ as
above, let us define a new order relation $\tri$ on
$\{0,1,\dots,n\}$. For any $p,q \in \{0,1,\dots,n\}$, $$ p \tri q,$$
if $s(p) < s(q)$ or if $s(p) = s(q)$ and $p>q$.

Additionally, let us define $m_i$ for $i = 0, 1, \dots, n$ such that
$m_i \in \{0,1,\dots,n\}$ and there are exactly $i$ elements $m\in
\{0,1,\dots,n\}$ such that $m \tri m_i$.

Now, let if we look at the $j$th cyclic shift of $\pi$, denoted
$\pi_j=a_{j+1}a_{j+2}\cdots a_na_1 \cdots a_j$, then the partial
sums of $\pi_j$, $s^j(p)$, is given by $s^j(p)= \sum_{i=1}^{p-j}
a_{i+j}=a_{i+1}+ a_{i+2} + \cdots + a_p$, where indices are
considered modulo $n$. This leads to the equation

\begin{equation}\label{E:sjp} s^j(p) = \begin{cases} s(p) - s(j), &
\text{ if $j \leq p \leq n$;}\\
s(p)-s(j)+k, & \text{ if $0 \leq p < j$.}\end{cases}\end{equation}

\begin{prop}
For any sequence $\pi$ with steps $a_i \in \{-1,1\}$ and sum $k=1$,
the $m_i$th cyclic shift of $\pi$ has exactly $i+1$ values of $p$
such that $s^{m_i}(p) \leq 0$.
\end{prop}

\begin{proof} Clearly, $s^{m_i}(m_i) = 0$, so there is at least one
such value of $p$. If $l < i$, then let us check $s^{m_i}(m_l)$. If
$s(m_l) = s(m_i)$ then $m_l > m_i$, so we know that $s^{m_i}(m_l) =
s(m_l)-s(m_i) = 0$. If, however, $s(m_l)<s(m_i)$ then $s^{m_i}(m_l)
\leq s(m_l) + 1 - s(m_i) < 1$, so $s^{m_i}(m_l) \leq 0$.

If $l>i$, then either $s(m_l) > s(m_i)$ or $m_l < m_i$ and $s(m_l) =
s(m_i)$. In the first case since $s^{m_i}(m_l) \geq s(m_l)-s(m_i)$,
clearly $s^{m_i}(m_l) > 0$. Otherwise, if $m_l < m_i$ then
$s^{m_i}(m_l) = 1 + s(m_l)-s(m_i)$ so since $s(m_l) = s(m_i)$,
$s^{m_i}(m_l) = 1$.

Thus there are exactly $i+1$ values of $p$, namely $m_0, m_1, \dots,
m_i$, such that $s^{m_i}(p) \leq 0$.
\end{proof}

The proof of the Cycle Lemma naturally follows from this
proposition.

\begin{proof}[Proof of Cycle Lemma]
Use the order relation $\tri$ to calculate $m_i$ for each $i = 0, 1,
\dots, n$ and let $\phi = \pi_{m_0}$. By definition, only
$s^{m_0}(m_0) = 0$. Since $s^{m_0}(p) = s_\phi(p-m_0)$, only
$s_\phi(0) \leq 0$, so all partial sums of $\phi$ for $p \geq 0$ are
positive.\end{proof}

Let $\pi$ represent a Dyck path from $(0,0)$ to $(2n+1,1)$, by
taking the path $(1,a_i)$. Let us restrict the ordering $\tri$ to
just the cyclic shifts of $\pi$ that begin with an up-step. Since
the order structure of $\{0,1, \dots, 2n+1\}$ is maintained we can
order the up-steps $j_k=m_i$ is the up-step with exactly $k$
up-steps $j$ such that $j \tri j_k$ for $k=0, 1, 2, \dots , n$.
Therefore, though there will be a total of $i$ vertices on or below
the $x$-axis in the $n_k$th cyclic shift of $\pi$, there will be
precisely $k$ up-steps that start on or below the $x$-axis.

Since there are $n+1$ total up-steps, there are $n+1$ cyclic shifts
of $\pi$ that start with an up-step. Now since there is precisely
one cyclic shift for every path, then we have that the total number
of paths with steps $a_i \in \{-1,1\}$ and sum $1$ that start with
an up-step and have exactly $k$ up-steps that start on or below the
$x$-axis is given by the fraction $\frac{1}{n+1}$ of the total
number of paths from $(0,0)$ to $(2n+1,1)$ that start with an
up-step. Thus it is the same as the number of paths from $(1,1)$ to
$(2n+1,1)$ which is precisely the number of paths from $(0,0)$ to
$(2n,0)$; that is, $\binom{2n}{n}$.

Therefore the number of paths from $(0,0)$ to $(2n+1,1)$ that start
with an up-step and have exactly $k$ up-steps that start on or below
the $x$-axis, for $k = 0, 1, \dots, n$, is given by
$\frac{1}{n+1}\binom{2n}{n}=C_n$. If we drop the initial up-step
then we are left with a path from $(0,0)$ to $(2n,0)$ with exactly
$k$ up-steps below the $x$-axis. So there are $C_n$ paths from
$(0,0)$ to $(2n,0)$ with exactly $k$ up-steps (and thus a total of
$2k$ steps) below the $x$-axis for each $k$.

\end{document}